\documentclass[11pt]{article}
\input epsf.tex
\usepackage{amssymb,latexsym,times}
\catcode`\@=11
\renewcommand\subsection{\@startsection{subsection}{2}{\z@}%
                                     {-3.25ex\@plus -1ex \@minus -.2ex}%
                                     {-0.01 mm}
                                     {\normalfont\large\bfseries}}

\catcode`\@=11
\renewcommand\subsubsection{\@startsection{subsubsection}{2}{\z@}%
                                     {-3.25ex\@plus -1ex \@minus -.2ex}%
                                     {-0.01 mm}
                                     {\normalfont\bfseries}}
\newtheorem{example}{Example}%[subsection]
\newtheorem{theorem}[example]{Theorem}%[section]
\newtheorem{corollary}[example]{Corollary}%[section]
\newtheorem{definition}[example]{Definition}%[section]
\newtheorem{proposition}[example]{Proposition}%[section]
\newtheorem{lemma}[example]{Lemma}%[section]
\def\resp{{\em resp.$\ $}}

\def\proof{\medskip\noindent {\it Proof --- \ }}

\def\cqfd{\hfill $\Box$ \bigskip}
\def\adots{\mathinner{\mkern2mu\raise1pt\hbox{.}
\mkern3mu\raise4pt\hbox{.}\mkern1mu\raise7pt\hbox{.}}}
\def\<{\langle\,}
\def\>{\,\rangle}

\def\ie{{\em i.e. }}

\def\sym{{\rm Sym}}

\def\SG{\mathfrak S}
\def\HH{\mathfrak H}

\def\l{\lambda}
\def\a{\alpha}
\def\de{\delta}
\def\b{\beta}
\def\ga{\gamma}
\def\N{{\mathbb N}}
\def\Z{{\mathbb Z}}
\def\C{{\mathbb C}}

\def\FF{{\mathbb F}}

\def\F{{\cal F}}
\def\P{{\cal P}}

\def\LL{{\mathfrak L}}

\def\SS{{\cal S}}

\def\mod{{\rm\ mod\ }}

\def\wt{{\rm wt}}
\def\g{\mathfrak g}

\def\Sl{\mathfrak{sl}}

\def\slchap{\widehat{\mathfrak{sl}}}
\def\glchap{\widehat{\mathfrak{gl}}}
\def\ASG{\widetilde{\mathfrak S}}

\def\A{{\cal A}}

\def\LL{{\cal L}}
\def\L{\Lambda}

\def\O{{\cal O}}
\def\bar{\overline}

\def\sp{{\rm spin}}
\def\<{\langle}
\def\>{\rangle}

\def\deg{{\rm deg}}

\def\CC{{\cal C}}
\def\m{\mu}

\def\le{\leqslant}
\def\ge{\geqslant}
\def\veps{\varepsilon}
\def\rad{{\rm rad}}
\def\Si{\Sigma}
\def\si{\sigma}
\def\th{s}
\def\vp{\varphi}
\def\inte{\rm int}
\def\ul{\underline{\l}}
\def\um{\underline{\m}}

\def\inf{\prec}
\def\sup{\succ}
\def\infe{\preceq}
\def\supe{\succeq}
\def\De{\Delta}
\newdimen\Squaresize \Squaresize=14pt
\newdimen\Thickness \Thickness=0.5pt
 
\def\Square#1{\hbox{\vrule width \Thickness
   \vbox to \Squaresize{\hrule height \Thickness\vss
      \hbox to \Squaresize{\hss#1\hss}
   \vss\hrule height\Thickness}
\unskip\vrule width \Thickness}
\kern-\Thickness}
 
\def\Vsquare#1{\vbox{\Square{$#1$}}\kern-\Thickness}

\title{\bf Some closed formulas for canonical bases \\ of Fock spaces}

\author{Bernard {\sc Leclerc} and Hyohe {\sc Miyachi}}

\date{April 2001}

\begin{document}
\maketitle

%%%%%%%%%%%%%%%%%%%%%%%%%%%%%%%%%%%%%%%%%%%%%%%%%%%%%%%%%%%%%%%%%%%
%  ABSTRACT
%%%%%%%%%%%%%%%%%%%%%%%%%%%%%%%%%%%%%%%%%%%%%%%%%%%%%%%%%%%%%%%%%%%
\vskip 1cm

\begin{abstract}\noindent
We give some closed formulas for certain vectors of the canonical bases
of the Fock space representation of $U_v(\slchap_n)$.
As a result, a combinatorial description of certain parabolic
Kazhdan-Lusztig polynomials for affine type $A$ is obtained.
\end{abstract}

\vskip 0.6cm
%%%%%%%%%%%%%%%%%%%%%%%%%%%%%%%%%%%%%%%%%%%%%%%%%%%%%%%%%%%%%%%%%%%
%  SECTION 1
%%%%%%%%%%%%%%%%%%%%%%%%%%%%%%%%%%%%%%%%%%%%%%%%%%%%%%%%%%%%%%%%%%%

\section{Introduction} \label{SECT1}
Let $\F_v$ be the Fock space representation of $U_v(\slchap_n)$
introduced by Hayashi \cite{H} and further studied by Misra and Miwa
\cite{MM}, Stern \cite{St}, and Kashiwara, Miwa and Stern \cite{KMS}. 
It has a standard basis $\Si=\{s(\l)\ | \ \l\in\P\}$ indexed by the set
$\P$ of all integer partitions.
In \cite{LT1} two canonical bases $B=\{G(\l)\ | \ \l\in\P\}$
and $B^-=\{G^-(\l)\ | \ \l\in\P\}$ of $\F_v$ have been constructed. 
The subset of $B$ consisting of the $G(\l)$'s for which $\l$
is $n$-regular coincides with Kashiwara's lower global basis
(or Lusztig's canonical basis)
of the irreducible sub-representation of $\F_v$ generated by the highest weight
vector $s(\emptyset)$.

The main motivation for introducing the bases $B$ and $B^-$
was their conjectural relation with the decomposition matrices of
the $q$-Schur algebras $\SS_m(q)$ defined by Dipper and James \cite{DJ}
in connection with the modular representation theory of the finite
groups $GL_m(\FF_q)$ in non-describing characteristic.
Conjecture 5.2 of \cite{LT1} was proved by Varagnolo and Vasserot
\cite{VV}, who established that the coefficients of the expansion
of $G^-(\l)$ on the basis $\Si$ are equal to the
Kazhdan-Lusztig polynomials appearing in Lusztig's character
formula for $U_\zeta(\Sl_r)$, where $\zeta$ is a complex primitive
$n$th root of 1.

Let $\ell$ be a prime number coprime to $q$ and such that the
multiplicative order of $q$ in $\overline{\FF}_\ell^*$ is equal to $n$.
To $w\in\N^*$ one associates a ``large $n$-core partition'' $\rho = \rho(w)$
(see below, Definition~\ref{CKcore}). 
Set $m=nw+|\rho|$, and let $B_{w,\rho}$ 
be the unipotent block of $\overline{\FF}_\ell GL_m(\FF_q)$ containing
the Specht-type modules $S(\l)$ labelled by partitions $\l$
with $n$-core $\rho$ and $n$-weight $w$ (see \cite{DJ}). 
Assume that $w < \ell$, and let 
$S(\l)=\rad^0(S(\l))\supset \rad^1(S(\l))\supset \rad^2(S(\l))\supset \ldots$ denote the 
radical series of $S(\l)$.  
In \cite{Mi} the graded composition multiplicities
\[
\rad_{\l,\m}(v) = \sum_{i\ge 0}\, [\rad^i(S(\l))/\rad^{i+1}(S(\l)) :
D(\mu)] \, v^i \]
of all $S(\l)$ in $B_{w,\rho}$ were computed explicitely in terms of
the Littlewood-Richardson coefficients,
using a Morita equivalence between $B_{w,\rho}$ and the principal
block of $\overline{\FF}_\ell GL_n(\FF_q)\wr \SG_w$ established
in \cite{HM}.
This Morita equivalence is similar to a Morita equivalence for
blocks of symmetric groups conjectured long ago by Rouquier 
(see \cite{R})
and recently proved by Chuang and Kessar \cite{CK}.
It was also conjectured in \cite{Mi} that
\[ 
\rad_{\l,\m}(v) = d_{\l,\m}(v)\,,
\]
where $d_{\l,\m}(v)$ denotes the $v$-decomposition number of
\cite{LT1}, that is, the coefficient of $s(\l)$ in the expansion
of $G(\m)$.
The main result of this paper is a proof of this conjecture.
Simultaneously, we also determine a similar expression for the coefficient
$e_{\l,\m}(-v^{-1})$ of $s(\m)$ in the expansion of~$G^-(\l)$.
As a consequence we obtain a closed and combinatorial expression for 
two families of parabolic affine Kazhdan-Lusztig polynomials.
It is remarkable that in fact all these polynomials are just monomials.  
Another consequence is that we obtain a representation-theoretical 
interpretation of these particular $v$-decomposition numbers :
the exponent of $v$ in $d_{\l,\m}(v)$ indicates to which layer of the radical
filtration of $S(\l)$ the copies of the simple module $D(\m)$ belong.

In terms of the Fock space, the block $B_{w,\rho}$ corresponds
to a distinguished weight space. More precisely, there is a
natural bijection between the set of $n$-core partitions and 
the set of extremal weights $\sigma(\L_0)$ of $\F_v$. 
Here $\L_0$ is the highest weight of $\F_v$ and $\sigma$ belongs to
the affine Weyl group $W=\ASG_n$ of $\slchap_n$. 
Let $\sigma_\rho(\L_0)$ be the extremal weight corresponding to $\rho$
in this bijection. Then our results give some closed formulas
in terms of the Littlewood-Richardson coefficients for the 
$\Si$-expansions of all the $G(\m)$'s and $G^-(\l)$'s of weight 
$\L_{w,\rho}:=\sigma_\rho(\L_0) - w\delta$, where $\delta$ is the imaginary root.   

Let $P(\F_v)$ denote the set of weights of $\F_v$.
We have the decomposition $P(\F_v) = \sqcup_{w\in\N} \O_w$, where
$\O_w = \{\L = \si\L_0 - w\de \ | \ \si\in W\}$ is the $W$-orbit of 
$\L_0-w\de$.
For $\L\in P(\F_v)$, let $T(\L)$ (\resp $T^-(\L)$) denote the transition
matrix from $\Si$ to $B$ (\resp from $\Si$ to $B^-$)
in the weight space $\F_v(\L)$. 
Having computed the matrices $T(\L_{w,\rho})$ and $T^-(\L_{w,\rho})$,
it is easy to determine $T(\L)$ and $T^-(\L)$ for many other weights $\L$.
Indeed, suppose that $\L\in P(\F_v)$ and $\alpha_i$ 
is a simple root of $\slchap_n$ such that 
$\L+\alpha_i  \not\in P(\F_v)$. 
Let $\si_i$ be the simple reflection of $W$
associated with $\alpha_i$.
Then $T(\L) = T(\si_i\L)$ and $T^-(\L) = T^-(\si_i\L)$.
It follows that, for each $w$, the orbit $\O_w$ 
can be splitted into a finite number of classes on which the
matrices $T(\L)$ and $T^-(\L)$ remain constant.
The class of $\L_{w,\rho}$ for which we have found closed formulas
is always infinite, hence for each $w$ our computations give the transition
matrices for an infinite number of weights of $\O_w$.
Moreover, in the case of $\slchap_2$ this class
is the only infinite one and therefore
our formulas calculate in this case the 
matrices for all but a finite number of
weights in each $\O_w$.  
These facts, which are easily deduced from the theory of crystal bases,
are the Fock space counterpart of Scopes' results \cite{Sc} about 
Morita equivalence for blocks of the symmetric groups, 
and of the analogues of these results for blocks of Hecke algebras and 
unipotent blocks of $GL_n(\FF_q)$ \cite{Jo1,Jo2}.

The paper is structured as follows.
In Section~\ref{SECT2} we review the correspondence between $n$-core
partitions and the $W$-orbit of the fundamental weight $\L_0$ of
$\slchap_n$.
Section~\ref{SECT3} and \ref{SECT4} recall the main facts about 
the Fock representations of $\slchap_n$ and $U_v(\slchap_n)$, and
introduce the canonical bases.
In Section~\ref{SECT5} we consider the space $\SS$ of symmetric
functions in $n$ independent sets of variables $A_0,\ldots ,A_{n-1}$
with coefficients in $\C(v)$.
The standard basis of $\SS$ is given by the products of Schur
functions $s_{\ul}=s_{\l^0}(A_0)\cdots s_{\l^{n-1}}(A_{n-1})$.
We introduce two new bases $\{\eta_{\ul}\}$ and
$\{\psi_{\ul}\}$ which are canonical with respect to
a certain bar-involution sending $v$ to $v^{-1}$ and two crystal 
lattices $\LL$ and $\LL^-$, and we calculate their expansions
on the standard basis in terms of the Littlewood-Richardson coefficients.
Heuristically, we regard $\SS$ (or rather its specialization $\SS_1$
at $v=1$) as the carrier space of the 
``canonical commutation relations'' representation of the homogeneous 
Heisenberg subalgebra of~$\glchap_n$. 
In fact, the sum of the homogeneous components of $\SS_1$ of degree
less than $\ell$ can be identified to the sum of the complexified
Grothendieck groups
\[
\bigoplus_{w<\ell} G(B(\overline{\FF}_\ell GL_n(\FF_q)\wr \SG_w))\,,
\]
where $B$ means the principal block. 
In this identification the vectors $s_{\ul}$ are mapped to 
the classes of the unipotent ordinary irreducible modules, 
and the vectors $\psi_{\ul}$ to the classes of the simple modules.
In Section~\ref{SECT6} we state our main result (Theorem~\ref{TH1}), namely
that the canonical bases of certain weight spaces of $\F_v$
coincide with the canonical bases of the corresponding 
homogeneous components of $\SS$ under a natural vector space
isomorphism. 
At $v=1$ this vector space isomorphism is essentially the natural intertwining
operator between the principal and the homogeneous realization of the
basic representation of $\slchap_n$ (see \cite{LL,Lei}). 
This is the counterpart in this setting of the Morita equivalence
of \cite{HM}, and it
gives immediately the above-mentioned formulas for $d_{\l,\m}(v)$
and $e_{\l,\m}(v)$
when $\l$ and $\mu$ have $n$-core $\rho$ and $n$-weight $\le w$
(Corollary~\ref{formule}).
Section~\ref{SECT7} is devoted to the proof of Theorem~\ref{TH1}.
Finally, Section~\ref{SECT8} develops in the context of
the Fock space the combinatorics underlying
Scopes' isometries between blocks of symmetric groups.

%%%%%%%%%%%%%%%%%%%%%%%%%%%%%%%%%%%%%%%%%%%%%%%%%%%%%%%%%%%%%
%    SECTION 2
%%%%%%%%%%%%%%%%%%%%%%%%%%%%%%%%%%%%%%%%%%%%%%%%%%%%%%%%%%%%%

\section{Combinatorics of partitions and the affine Weyl group} \label{SECT2}

\subsection{}\label{abacus} 
A partition is a finite non-increasing sequence of positive integers.
We shall denote by $\P$ the set of all partitions.
By convention $\P$ contains the empty partition $\emptyset$.
To a partition 
\[
\lambda = (\lambda_1, \ldots , \lambda_k)
\]
one associates the infinite decreasing sequence of $\beta$-numbers
\[
\beta(\lambda) = 
(\lambda_1, \lambda_2 -1, \lambda_3 - 2, \ldots , \lambda_i-i+1, \ldots )\,,
\]
where $\lambda_i$ is assumed to be $0$ for $i>k$.

\begin{figure}[t]
\begin{center}
\leavevmode
\epsfxsize =10cm
\epsffile{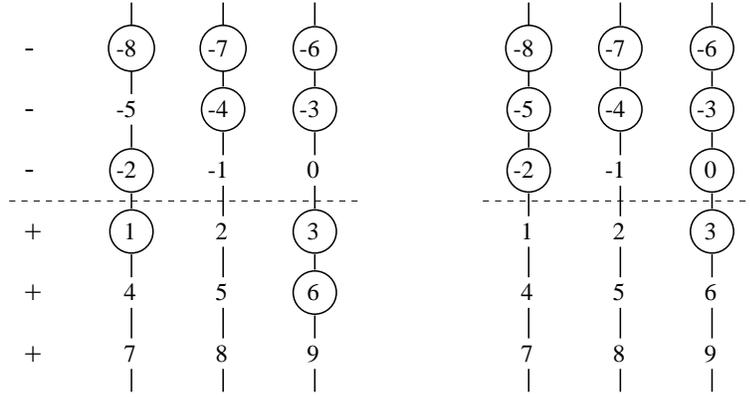}
\end{center}
\caption{\label{Fig1} \small The $3$-abacus representation of
$\lambda = (6,4,3,1,1,1)$ and of its $3$-core 
$\lambda_{(3)} = (3,1)$.}
\end{figure}

Fix an integer $n\ge 2$, and form an infinite abacus with $n$ 
runners labelled $1, \ldots , n$ from left to right.
The positions on the $i$th runner are labelled by the 
integers having residue $i$ modulo $n$.
By placing a bead on each $\beta$-number of $\lambda$ one 
gets the abacus representation of~$\lambda$.
As is well known, sliding all the beads in the $n$-abacus
representation of $\lambda$ as high as they will go
produces the $n$-core $\lambda_{(n)}$ of $\lambda$.
This is illustrated in Figure~\ref{Fig1} for 
$\lambda = (6,4,3,1,1,1)$ and $n=3$.
In that case $\beta(\lambda) = (6, 3, 1, -2, -3, -4, -6, -7, -8,
\ldots)$.

Note that in the abacus representation of a partition,
the number of occupied positive positions is always equal
to the number of vacant non-positive positions.
In particular, we see that the $n$-core partitions are
in one-to-one correspondence with the $n$-tuples of
integers $(a_1,\ldots , a_n)$ such that $\sum_i a_i = 0$.
Thus, as shown in Figure~\ref{Fig1}, to the $3$-core $(3,1)$
is associated the $3$-tuple $(0,-1,1)$.

The next definition introduces certain ``large $n$-core partitions''.
They were first considered by Rouquier as labels of certain
good blocks of the symmetric groups.

\begin{definition}\label{CKcore}{\rm
Let $w\in\N^*$. The partition $\rho = \rho(w)$ associated
with $w$ is the $n$-core corresponding to the $n$-tuple 
\[
(a_1,\ldots , a_n) = \left\{
\matrix{
\left(\frac{(1-n)(w-1)}{2},\frac{(3-n)(w-1)}{2},\ldots , 
\frac{(n-3)(w-1)}{2},\frac{(n-1)(w-1)}{2}\right)
\mbox{\ if $n$ or $w$ is odd,} \cr
\left(\frac{w}{2},\frac{3w-2}{2},\ldots,\frac{(n-1)w-n+2}{2},
\frac{(1-n)w+n-2}{2},\ldots,\frac{2-3w}{2},-\frac{w}{2}\right)
\mbox{\ otherwise.\ \ \ }
}\right.
\] 
}
\end{definition}
\begin{example}
{\rm For $n=4$ and $w=3$, we have 
\[
(a_1,a_2,a_3,a_4)=(-3,-1,1,3)
\quad\mbox{and}\quad 
\rho = (12,9,6^2,4^2,2^3,1^3).
\]
For $n=4$ and $w=4$, we have 
\[
(a_1,a_2,a_3,a_4)=(2,5,-5,-2)
\quad\mbox{and}\quad
\rho = (18,15,12,9^2,7^2,5^2,3^3,2^3,1^3).
\]
}
\end{example}

\subsection{} Let $\g = \slchap_n$ be the affine Lie algebra of type 
$A_{n-1}^{(1)}$ \cite{Ka}.
Denote by $\Lambda_0, \ldots , \Lambda_{n-1}$ the fundamental weights
of $\g$, by $\alpha_0, \ldots , \alpha_{n-1}$ its simple roots, and
by $\delta$ the imaginary root. They are related by
\[
\left\{
\begin{array}{llll}
\alpha_0 &=& 2\Lambda_0 - \Lambda_1 - \Lambda_{n-1} + \delta,
  \cr
\alpha_i &=& -\Lambda_{i-1} + 2\Lambda_i -\Lambda_{i+1} &
(i=1,\ldots,n-2), \cr
\alpha_{n-1} &=& -\Lambda_0 -\Lambda_{n-2} + 2\Lambda_{n-1}.
\end{array}
\right.
\]
if $n\ge 3$, and
\[
\left\{
\begin{array}{lll}
\alpha_0 &=& 2\Lambda_0 - 2\Lambda_1  + \delta,  \\
\alpha_1 &=& -2\Lambda_0 + 2\Lambda_1.
\end{array}
\right.
\]
if $n=2$.
The free $\Z$-modules  
$P = \left(\oplus_{i=0}^{n-1}\Z \Lambda_i\right)\oplus \Z\delta$  
and $Q = \oplus_{i=0}^{n-1}\Z \alpha_i$ are called respectively
the weight lattice and the root lattice of $\g$.
Let $a_{ij}$ be the coefficient of $\Lambda_j$ in $\alpha_i$.
The $(n-1)\times (n-1)$ matrix $A=(a_{ij})$ is called the
Cartan matrix of $\g$.
One defines a non-degenerate symmetric bilinear form on $P$ by
\[
\left\{
\begin{array}{llll}
(\alpha_i , \alpha_j) &=& a_{ij} & (0\le i,j \le n-1), \cr
(\Lambda_0,\Lambda_0)&=&0,\cr
(\Lambda_0,\alpha_0) &=& 1, \cr
(\Lambda_0 , \alpha_i) &=& 0 & (1 \le i \le n-1).
\end{array}
\right.
\]
The Weyl group $W$ of $\g$ is the subgroup of $GL(P)$ generated by
the simple reflections $\si_i$ defined by
\[
\si_i(\Lambda) = \Lambda - (\alpha_i , \Lambda) \alpha_i
\qquad (\Lambda \in P,\ 0\le i \le n-1).
\]
This is a Coxeter group of type $A_{n-1}^{(1)}$, and it is isomorphic
to the semi-direct product of the symmetric group $\SG_n$
by the root lattice $Q_0 = \oplus_{i=1}^{n-1} \Z \alpha_i$ of the
finite-dimensional Lie algebra $\g_0 = \Sl_n$. 

\subsection{} \label{X}
It will be convenient to realize $P$ and $Q$ as sublattices of a
lattice $X$ of rank $n+2$. Namely we set 
\[
X = \left(\oplus_{i=1}^{n}\Z \veps_i\right)\oplus \Z\Lambda_0 \oplus
\Z\delta
\]
and we identify the simple roots and fundamental weights to the
following elements of $X$:
\[
\alpha_i = \veps_i - \veps_{i+1} \quad (i=1,\ldots ,n-1), \qquad
\alpha_0 = \veps_n - \veps_1 + \delta, 
\]
\[
\Lambda_i = \veps_1 + \cdots + \veps_i - 
\frac{i}{n}(\veps_1+ \cdots + \veps_n) + \Lambda_0
\qquad (i=1,\ldots n-1).
\]
Note that the root lattice $Q_0$ of $\g_0$ gets identified in this realization
to 
\[
Q_0 = \left\{ \sum_{i=1}^n a_i \veps_i \in X \ ; \ \sum_{i=1}^n a_i = 0 \right\}.
\]  
It is easy to check that the above bilinear form is the
restriction to $P$ of the bilinear form on $X$ given by
\[
\left\{
\begin{array}{llll}
(\veps_i,\veps_j) &=& \delta_{ij} & (1\le i,j \le n), \cr
(\Lambda_0,\Lambda_0) &=& 0, \cr
(\Lambda_0,\veps_i) &=& 0 & (1\le i \le n), \cr
(\delta,\delta) &=& 0, \cr
(\delta,\veps_i) &=& 0 & (1\le i \le n), \cr
(\Lambda_0,\delta) &=& 1.
\end{array}
\right.
\]
This implies that the action of $W$ on $P$ can be lifted to $X$
by setting
\[
\left\{
\begin{array}{llll}
\si_0(\Lambda_0) &=& \Lambda_0 - \alpha_0, \cr
\si_0(\delta) &=& \delta, \cr
\si_0(\veps_1) &=& \veps_n + \delta, \cr
\si_0(\veps_n) &=& \veps_1 - \delta, \cr
\si_0(\veps_j) &=& \veps_j, & (j \not =1,n)
\end{array}
\right.
\quad
\left\{
\begin{array}{llll}
\si_i(\Lambda_0) &=& \Lambda_0, \cr
\si_i(\delta) &=& \delta, \cr
\si_i(\veps_i) &=& \veps_{i+1}, \cr
\si_i(\veps_{i+1}) &=& \veps_i, \cr
\si_i(\veps_j) &=& \veps_j, & (j \not = i, i+1)
\end{array}
\right.
\quad
(i \not = 0).
\]

\subsection{} \label{n-cores}
The stabilizer $W_0$ of $\Lambda_0$ is the subgroup of $W$ generated by
$\si_1, \ldots , \si_{n-1}$.
Therefore, it is isomorphic to the symmetric group $\SG_n$.
Hence, the orbit $W\Lambda_0 \simeq W/W_0$ is in a natural one-to-one
correspondence with the root lattice $Q_0$ of $\g_0$.
Namely, for $\sigma \in W$ write
\[
\sigma(\Lambda_0) = \Lambda_0 + d\delta + \sum_{i=1}^n a_i \veps_i\,.
\]
By the formulas of \ref{X} for the action of $W$ on $X$, 
we see that the integer coordinates $a_i$ satisfy
$\sum_i a_i = 0$, and that the projection 
$\sigma(\Lambda_0) \mapsto \sum_{i=1}^n a_i \veps_i$ induces a bijection
from $W\Lambda_0$ to $Q_0$.

\begin{figure}[t]
\begin{center}
\leavevmode
\epsfxsize =10cm
\epsffile{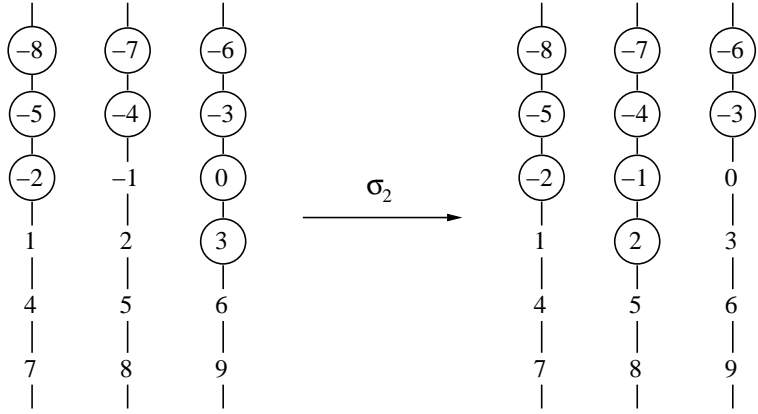}
\end{center}
\caption{\label{Fig2} \small The $3$-core partition $(3,1)$ is mapped under
$\si_2$ to the $3$-core partition $(2)$.}
\end{figure}

Comparing with \ref{abacus}, 
we obtain that the orbit $W\Lambda_0$
can be identified with the set $\CC_n$ of $n$-core partitions. 
In other words, there is a natural action of the affine Weyl group
$W$ on~$\CC_n$, which can be easily described using the abacus
representation of the elements of~$\CC_n$.
Namely, since $W_0$ acts on $Q_0$ by permuting the vectors $\veps_i$,
we see that for $i=1,\ldots , n-1$, the action of $\si_i$ on an
$n$-core partition is obtained by switching the beads of the $i$th
and $(i+1)$th runners of its abacus representation.
This is illustrated in Figure~\ref{Fig2}.
On the other hand, we have
\[
\si_0\left( \Lambda_0 + d\delta + \sum_{i=1}^n a_i \veps_i \right)
= \Lambda_0 +(d-1)\delta + (a_n+1)\veps_1 + a_2\veps_2 + \cdots + 
a_{n-1}\veps_{n-1} + (a_1-1)\veps_n\,,
\]
so that the action of $\si_0$ on an
$n$-core partition is obtained by switching the beads of the first
and last runners of its abacus representation,  adding one more
bead to the first runner and removing one bead from the last runner.
%
%%%%%%%%%%%%%%%%%%%%%%%%%%%%%%%%%%%%%%%%%%%%%%%%%%%%%%%%%%%%%
%    SECTION 3
%%%%%%%%%%%%%%%%%%%%%%%%%%%%%%%%%%%%%%%%%%%%%%%%%%%%%%%%%%%%%

\section{The Fock space representation of $\g$} 
\label{SECT3}

\subsection{} 
Let $V(\Lambda_0)$ denote the irreducible $\g$-module
with highest weight $\Lambda_0$ (\cite{Ka}, 9). 
It is the simplest infinite-dimensional $\g$-module, 
and as such, it is often called the
basic representation of~$\g$.
The easiest way to construct it is to embed it into a larger
representation of $\g$ called the Fock space representation
\cite{JiMi}.
The Fock space $\F$ has a distinguished $\C$-basis 
$\Si=\{ s(\lambda)\ ;\ \lambda \in \P\}$, on which the 
Chevalley generators $e_i, f_i \ (0\le i \le n-1)$ of $\g$
act by
\begin{equation}\label{Robinson}
e_i s(\lambda) = \sum_\mu s(\mu),
\qquad
f_i s(\lambda) = \sum_\nu s(\nu),
\end{equation}
where the first sum is over all partitions $\mu$ obtained
by removing from the Young diagram of $\lambda$ a node with
$n$-residue $i$, and the second sum is over all $\nu$'s
obtained from $\lambda$ by adding an $i$-node.     
Denoting by $\emptyset$ the empty partition, we see immediately
that $s(\emptyset)$ is a highest weight vector of $\F$ with
weight $\Lambda_0$. The cyclic submodule $U(\g)\,s(\emptyset)$
is isomorphic to $V(\Lambda_0)$.

These infinite-dimensional 
$\C$-vector spaces are graded by putting 
$\deg\, s(\lambda) = m$ if $\lambda$
is a partition of $m$. 
The dimension of the degree $m$ component of $V(\Lambda_0)$ is known to be
equal to the number of $n$-regular partitions of $m$.

For a $\g$-module $M$ and a weight $\Lambda\in P$, we denote
by $M(\Lambda)$ the $\Lambda$-weight space of $M$.
Let $P(M)$ denote the set of weights of $M$,
that is, the set of $\Lambda \in P$ such that $M(\Lambda)$ 
is non-zero.
It is known (see \cite{Ka}, 12.6) that 
\[
P(V(\Lambda_0)) = P(\F) = \{ \sigma(\Lambda_0) - w\delta \ ; \ 
\sigma \in W,\ w \in \N \}.
\]
Recall from \ref{n-cores} that there is a natural bijection,
say $\gamma$, from the orbit $W\Lambda_0$ to the set $\CC_n$ of $n$-cores.
Then the weight-space $\F(\L)$ associated with the
weight $\L = \si(\L_0) - w\delta$ consists precisely
of the span of the vectors $s(\l)$ where $\l$
has $n$-core $\gamma(\si(\L_0))$
and $n$-weight $w$.

\subsection{} \label{integrable} 
The Fock space is an integrable $\g$-module,
that is, the operators $e_i$ and $f_i$ are locally nilpotent on $\F$.
This implies that 
$\dim \F(\Lambda) = \dim \F({\sigma(\Lambda)})$
for all $\Lambda \in P(\F)$ and $\sigma \in W$.
In fact, it is well known that 
\[
\dim \F(\si(\L_0)-w\delta) =
\sum_{k_1+\cdots +k_n=w}
p(k_1)\cdots p(k_n)
\]
where $p(k)$ is the number of integer partitions of $k$.

Define the following linear operators on $\F$ :
\[
r_i = \exp(e_i)\,\exp(-f_i)\,\exp(e_i) 
\qquad (0\le i \le n-1),
\] 
(see \cite{Ka} 3.8).
They satisfy
$r_ie_i=-f_ir_i$ and $r_if_i=-e_ir_i$, hence
\[
r_i^{-1} = \exp(f_i)\,\exp(-e_i)\,\exp(f_i).
\]
Moreover, 
$r_i(\F(\Lambda)) = \F(\sigma_i(\Lambda))$
for all $\Lambda \in P(\F)$. 

It is easy to see that the divided powers $e_i^k/k!$ and $f_i^k/k!$
preserve the $\Z$-lattice $\F_\Z$ spanned by the basis
$\{ s(\lambda)\ ;\ \lambda \in \P\}$.
It follows that $r_i$ also preserves $\F_\Z$ for all $i$.
 
\subsection{} \label{scalar}
$\F$ is endowed with a canonical scalar product
defined by 
\[
\< s(\lambda)\,,\,s(\mu)\> = \delta_{\lambda,\mu}
\qquad (\lambda, \mu \in \P).
\]
This form is compatible with the actions of $\g$ and of the $r_i$, in the sense that
\[
\<e_i x\,,\,y\> = \<x \,,\,f_i y\>, 
\quad
\<r_i x\,,\,y\> = \<x \,,\,r_i^{-1} y\> 
\qquad (x,y \in \F,\ 0\le i \le n-1).
\]
Hence we see that $r_i$ is an orthogonal transformation of $\F$.

\subsection{}\label{bosons}
For $k\in \N^*$, define an endomorphism $d_k$ of $\F$ by
\[
d_k\, s(\l) = \sum_\m (-1)^{\sp(\m/\l)}\,s(\m)
\]
where the sum is over all $\m$ such that $\m/\l$ is a horizontal
$n$-ribbon strip of weight $k$, and $\sp(\m/\l)$ is the spin of the
strip (see \cite{LT2}).
It is known that the $d_k$ pairwise commute, and also commute with
the action of $\g' = [\g,\g]$ on $\F$.
Moreover $\F$ becomes a cyclic module under the action
of $\A := U(\g')\otimes\C[d_k\,;\,k\ge 1]$, that is,  
$\F \simeq \A \, s(\emptyset)$.
    
%%%%%%%%%%%%%%%%%%%%%%%%%%%%%%%%%%%%%%%%%%%%%%%%%%%%%%%%%%%%%%
%    SECTION 4
%%%%%%%%%%%%%%%%%%%%%%%%%%%%%%%%%%%%%%%%%%%%%%%%%%%%%%%%%%%%%

\section{The Fock space representation of $U_v$} \label{SECT4} 

\subsection{} Let $U_v = U_v(\slchap_n)$ be the quantum affine algebra
of type $A_{n-1}^{(1)}$ over the field $K:=\C(v)$. 
Let $\F_v$ be the $v$-deformation of the Fock space introduced
by Hayashi \cite{H} and further studied by Misra and Miwa \cite{MM}.
As a vector space 
$\F_v = K\otimes_\C \F$, and $\Si=\{s(\l)\,|\,\l\in\P\}$
is a $K$-basis of $\F_v$.
The Chevalley generators $E_i, F_i \ (0\le i \le n-1)$
of $U_v$ act on $\F_v$ via some simple $v$-analogues 
of the formulas (\ref{Robinson}), as follows.

Let $\l$ and $\m$ be two Young diagrams
such that $\m$ is obtained from $\l$ by adding an $i$-node $\gamma$.
Such a node is called a removable $i$-node of $\m$, or
an indent $i$-node of $\l$.
Let $I_i^r(\l,\m)$ (\resp $R_i^r(\l,\m)$)
be the number of indent $i$-nodes of $\l$
(\resp of removable $i$-nodes of $\l$)
situated to the right of $\gamma$ ($\gamma$ not included).
Set
$N_i^r(\l,\m)=I_i^r(\l,\m)-R_i^r(\l,\m)$. Then
\begin{equation}\label{ACTF}
F_i s(\l) = \sum_\m v^{N_i^r(\l,\m)} s(\m),
\end{equation}
where the sum is over all partitions $\m$ such that
$\m/\l$ is an $i$-node.  
Similarly
\begin{equation}\label{ACTE}
E_i s(\m) = \sum_\l v^{-N_i^l(\l,\m)} s(\l),
\end{equation}
where the sum is over all partitions $\l$ such that
$\m/\l$ is an $i$-node, and $N_i^l(\l,\m)$ is defined
as $N_i^r(\l,\m)$ but replacing right by left.  

\subsection{} \label{q-bosons}
The cyclic submodule $U_v\,s(\emptyset)$ is isomorphic to the
basic representation $V_v(\Lambda_0)$ of $U_v$.
As in \ref{bosons}, define operators $D_k\ (k\ge 1)$ acting on $\F$ by
\[
D_k \, s(\l) = \sum_\m (-v)^{-\sp(\m/\l)} \, s(\m)\,,
\]
the sum being over all $\m$ such that $\m/\l$ is a horizontal
$n$-ribbon strip of weight $k$.
They pairwise commute and they also commute with the action
of the subalgebra $U'_v$ of $U_v$ obtained by omitting the
degree generator of the Cartan part.
The Fock space $\F_v$ becomes a cyclic module 
under the action of 
$\A_v := U'_v \otimes K[D_k\,;\,k\ge 1]$,
generated by $s(\emptyset)$.

\subsection{} \label{SECT4.3}
By Kashiwara's theory of crystal bases \cite{Ka1,Ka2,Ka3},
$V_v(\Lambda_0)$ has two canonical bases dual to each other
called the lower global base and the upper global base. 

In \cite{LT1} (see also \cite{LT2}) the canonical bases 
$B=\{G(\lambda)\}$ and 
$B^-=\{G^-(\lambda)\}$ of $\F_v$
were introduced.
We shall review their definition.

Let $x \mapsto \overline{x}$ be the bar involution
of $\F_v$ defined in \cite{LT1}. 
It is the unique semi-linear map satisfying 
$\overline{s(\emptyset)} = s(\emptyset)$ and
\[
\overline{F_i\,x} = F_i \overline{x}\,,
\qquad
\overline{D_k\,x} = D_k \overline{x}\,,
\qquad
(0\le i \le n-1,\ k\ge 1,\ x\in\F_v).
\]
Let $L$ (\resp $L^-$) be the free $\Z[v]$-module 
(\resp $\Z[v^{-1}]$-module) with basis ${s(\lambda)}$.
The bases $B$ and $B^-$ are characterized by the 
properties
\[
\overline{G(\lambda)} = G(\lambda), 
\quad
\overline{G^-(\lambda)} = G^-(\lambda),
\quad
G(\lambda) \equiv s(\lambda) \mbox{\ \mod\ } vL,
\quad
G^-(\lambda) \equiv s(\lambda) \mbox{\ \mod\ } v^{-1}L^-,
\]
for $\lambda\in\P$.
The subset of $B$ consisting of the $G(\l)$ for which $\l$
is $n$-regular coincides with the lower global base of $V_v(\L_0)$.

\subsection{}
Define polynomials $d_{\lambda,\mu}(v)$ and $e_{\l,\m}(v)$ by
\[
G(\mu) = \sum_\lambda d_{\lambda,\mu}(v) s(\lambda)
\]
and
\[
G^-(\l) = \sum_\m e_{\l,\m}(-v^{-1}) s(\m).
\]
The polynomials $d_{\lambda,\mu}(v)$ and $e_{\l,\m}(v)$ are known
to be parabolic Kazhdan-Lusztig polynomials and they belong to
$\N[v]$ (see \cite{VV,LT2,KT}).

%%%%%%%%%%%%%%%%%%%%%%%%%%%%%%%%%%%%%%%%%%%%%%%%%%%%%%%%%%%%%
% SECTION 5
%%%%%%%%%%%%%%%%%%%%%%%%%%%%%%%%%%%%%%%%%%%%%%%%%%%%%%%%%%%%%
%
\section{Symmetric functions and Heisenberg algebras} \label{SECT5}
\subsection{}\label{SECT5.1}
Let $\sym$ denote the $K$-algebra of symmetric functions in
a countable set $X$ of indeterminates  (see \cite{Mcd}).
It is known that $\sym$ is a polynomial algebra in the 
power-sums symmetric functions 
\[
p_k = \sum_{x\in X} x^k \,, \qquad (k \ge 1).
\]
Other systems of generators are the elementary symmetric functions
$e_k$ and the complete symmetric functions $h_k$. 
The Schur functions $s_\lambda$ form a linear basis of $\sym$,
and we denote by $\< \cdot\, , \cdot \>$ the scalar
product for which this basis is orthonormal.
We may identify $\sym$ to the Fock space representation $\F_v$ of
$U_v(\slchap_n)$
by identifying the basis $\{s_\lambda\}$ with the standard basis $\Si$
of $\F_v$.

For $f\in\sym$, denote by $\hat{f}$ the operator of multiplication 
by $f$ and by $D_f$ the adjoint operator, that is,
\[
\<D_f(g) \,,\, h\> = \<g\,,\,f h\> \qquad (f,g,h \in \sym).
\]
It is known that $D_{p_k} = k\,(\partial/\partial p_k)$ 
(see \cite{Mcd}, I.5 Ex.3). 
It follows that the operators $\hat{f}$, $D_f$ $(f\in\sym)$ generate
the enveloping algebra of the Heisenberg Lie algebra $\HH$ spanned by
the operators $\hat{p}_k$, $\partial/\partial p_k$ $(k\in\N^*)$
and the identity.
Thus, $\sym \simeq \F_v$ can be regarded as the
``canonical commutation relations'' representation of $\HH$.
 
\subsection{}\label{SECT5.2}
We introduce the generating function
\[
H(t,X) = \sum_{k\ge0} h_k t^k = \prod_{x\in X} (1-xt)^{-1}
\]
and we put $H(X) := H(1,X)$.
Then, if we denote by $Y$ a second countable set of indeterminates,
and write $XY=\{xy \ |\ x\in X, y\in Y\}$, we have the Cauchy identity
(see \cite{Mcd} I.4.3)
\begin{equation}\label{CAUCHY}
\sum_{\l\in\P} s_\l(X)s_\l(Y) = H(XY) = \prod_{x\in X,y\in Y} (1-xy)^{-1}.
\end{equation}
More generally, if $U=\{u_\l\,,\,\l\in\P\}$, $V=\{v_\l\,,\,\l\in\P\}$
are two homogeneous bases of $\sym$ with 
$\deg\, u_\l = \deg\, v_\l = |\l|$,
then $U$ and $V$ are adjoint to each other with
respect to $\< \cdot\, , \cdot \>$ if and only if
\[
\sum_{\l\in\P} u_\l(X)v_\l(Y) = H(XY)
\]   
(see \cite{Mcd} I.4.6).

\subsection{}\label{SECT5.3}
Let $A_0, \ldots , A_{n-1}$ be $n$ countable sets of
indeterminates.
We denote by 
\[
\SS = \sym(A_0,\ldots ,A_{n-1})
\]
the algebra 
over $K$ of functions symmetric in each 
set $A_0, \ldots ,A_{n-1}$ separately.
This is the polynomial algebra in the variables
$p_k(A_i)\ (k\in\N^*, \ 0\le i \le n-1)$.
A linear basis of $\SS$ is given by the products
\[
s_{\ul} := s_{\l^{0}}(A_0)\cdots s_{\l^{n-1}}(A_{n-1}),
\qquad
\ul=(\l^{0},\cdots ,\l^{n-1})\in\P^n.
\]
The space $\SS$ carries a scalar product defined
by
\[
\< f_0(A_0)\cdots f_{n-1}(A_{n-1})\,,\, g_0(A_0)\cdots g_{n-1}(A_{n-1})\>
=\< f_0,g_0\>\cdots \< f_{n-1},g_{n-1}\>
\]
for $f_0,\ldots ,f_{n-1},g_0,\ldots ,g_{n-1} \in \sym$.
Clearly, the basis $\{s_{\ul}\ |\ \ul \in \P^n\}$ is orthonormal
for this scalar product.

As in \ref{SECT5.1}, we regard $\SS$ as the ``canonical commutation
relations'' representation of the Heisenberg algebra $\HH_n$
generated by the operators
\[
\hat{p}_k(A_i), \quad \frac{\partial}{\partial p_k(A_i)}, \qquad 
(k\in\N^*, \ 0\le i \le n-1).
\]

In the rest of this section we are going to introduce two
canonical bases of $\SS$.

\subsection{}
Given $q_1(v),\ldots ,q_k(v) \in \Z[v,v^{-1}]$ and
$i_1,\ldots ,i_k \in \{0,\ldots ,n-1\}$, we define symmetric functions
of the formal set of variables $q_1(v)A_{i_1} + \cdots + q_k(v)A_{i_k}$  
as follows.
For $f\in \sym$  
we let $f(q_1(v)A_{i_1} + \cdots + q_k(v)A_{i_k})$ 
denote the image of $f$ by
the algebra homomorphism from $\sym$ to $\SS$ which maps
$p_k$ to 
\[
p_k(q_1(v)A_{i_1} + \cdots +q_k(v)A_{i_k}):=
q_1(v^k)p_k(A_{i_1}) + \cdots +q_k(v^k)p_k(A_{i_k})\,.
\]
(For the reader familiar with the language of $\l$-rings \cite{Kn}, we
consider $A_0,\ldots , A_{n-1},v$ as elements of a
$\l$-ring, $v$ being invertible of rank 1.)
For example, for $j\in \{0,\ldots ,n-2\}$ we define
\[
p_k(A_j+vA_{j+1}+\cdots +v^{n-1-j}A_{n-1}):=
p_k(A_j)+v^kp_k(A_{j+1}) + \cdots + v^{(n-j-1)k}p_k(A_{n-1})\,.
\]    
Then we have (see \cite{Mcd} I.5.9)
%\begin{equation}
%\begin{multiline}
\[
s_\lambda(A_j+vA_{j+1}+\cdots +v^{n-1-j}A_{n-1}) 
\qquad\qquad\qquad\qquad\qquad\qquad\qquad\qquad\qquad
\]
\[
\\
\qquad
= 
\sum_{\alpha^{j},\ldots ,\alpha^{n-1}\in\P} 
v^{\sum_{j\le i\le n-1} (i-j)|\alpha^{i}|}\, 
c_{\alpha^{j},\ldots ,\alpha^{n-1}}^\lambda
 \,s_{\alpha^{j}}(A_j)\cdots s_{\alpha^{n-1}}(A_{n-1})\,,
\]
%\end{multiline}
%\end{equation}
where the $c_{\alpha^{j},\ldots ,\alpha^{n-1}}^\lambda$ are the 
Littlewood-Richardson coefficients, that is,
\begin{equation}
c_{\alpha^{j},\ldots ,\alpha^{n-1}}^\lambda
:= \<s_{\a^j}\cdots s_{\a^{n-1}}\,,\, s_\l\>\,.
\end{equation}
Similarly, for $j=1,\ldots ,n-1$, we set
\[
p_k(A_j-vA_{j-1}):=p_k(A_j)-v^kp_k(A_{j-1})\,,
\]
and we have (see \cite{Mcd} I.3.10) 
\begin{equation}\label{EQ3}
s_\lambda(A_j-vA_{j-1})
= \sum_{\b\in\P} s_{\l/\b}(A_j) s_\b(-vA_{j-1})  
= \sum_{\alpha,\beta\in\P}(-v)^{|\beta|}\, c_{\alpha,\beta}^\lambda
 \,s_\alpha(A_j) s_{\beta'}(A_{j-1})\,,
\end{equation}
where 
$\beta'$ denotes the partition conjugate to $\beta$.

\subsection{} For $\ul = (\l^{0},\ldots ,\l^{n-1})\in\P^n$
we write $|\ul|=\sum_{0\le i\le n-1} |\l^{i}|$.
Consider two bases of $\SS$
$U=\{u_{\ul}\,,\,\ul\in\P^n\}$ and $V=\{v_{\ul}\,,\,\ul\in\P^n\}$
consisting of homogeneous elements of degree 
$\deg\, u_{\ul} = \deg\, v_{\ul} = |\ul|$.
Then \ref{SECT5.2} may clearly be generalized, as follows.
\begin{lemma}\label{LEM3}
The bases $U$ and $V$ are adjoint to each other with
respect to $\< \cdot\, , \cdot \>$ if and only if
\[
\sum_{\ul\in\P^n} 
u_{\ul}(A_0,\ldots,A_{n-1})\,v_{\ul}(B_0,\ldots,B_{n-1}) 
= H(A_0B_0)\cdots H(A_{n-1}B_{n-1}).
\] 
\end{lemma}

\subsection{} For $\ul = (\l^0,\ldots,\l^{n-1}) \in \P^n$, 
we define the following elements of $\SS$:
\begin{eqnarray*}
\eta_{\ul}(v) &=&
s_{\l^{0}}(A_0)\,s_{\l^{1}}(A_1-vA_0)\,s_{\l^{2}}(A_2-vA_1)
\ \cdots\  
s_{\l^{n-1}}(A_{n-1}-vA_{n-2})\,,
\\  
\varphi_{\ul}(v) &=& s_{\l^{0}}(A_0+vA_1+\cdots +v^{n-1}A_{n-1})
\,s_{\l^{1}}(A_1+vA_2+\cdots +v^{n-2}A_{n-1})
\\
&&\cdots\ 
s_{\l^{n-3}}(A_{n-3}+vA_{n-2}+v^2A_{n-1})\,
s_{\l^{n-2}}(A_{n-2}+vA_{n-1})\,
s_{\l^{n-1}}(A_{n-1}).
\end{eqnarray*}
It is easy to see that these are two bases of $\SS$.
By \ref{SECT5.3}, their expansions on the orthonormal basis
$\{s_{\ul}\}$ are readily computed
in terms of the Littlewood-Richardson coefficients.
In particular, 
\begin{lemma} \label{LEM4} For $\ul,\um \in \P^n$ we have
\[
\< s_{\ul} \,,\,\eta_{\um}(v) \> =
(-v)^{\de(\ul,\um)}
\sum_{\a^{0},\ldots,\a^{n} \atop \b^{0},\ldots,\b^{n-1}}
%c_{\m^{0}(\a^{1})'}^{\l^{0}}
%c_{\a^{n-1}\l^{n-1}}^{\m^{n-1}}
\prod_{0\le j \le n-1}
c_{\a^{j}\b^{j}}^{\m^{j}}
c_{\b^{j}(\a^{j+1})'}^{\l^{j}}\,,
\]
where $\a^{0},\ldots,\a^{n},\b^{0},\ldots,\b^{n-1}$ run through
$\P$ subject to the conditions
\[
|\a^i| = \sum_{0\le j \le i-1} |\l^j|-|\m^j|\,, 
%\quad (0\le i \le n)
\qquad
|\b^i| = |\m^i| + \sum_{0\le j \le i-1} |\m^j|-|\l^j|\,,
%\quad (0\le i \le n-1),
\] 
and 
\[
\de(\ul,\um) := \sum_{0\le j \le n-2}(n-1-j)(|\l^j|-|\m^j|)\,.
\]
\end{lemma}
Here we have used the convention that an empty sum
is equal to 0. Thus $|\a^0|=0$ and $|\b^0| = |\m^0|$, so that
$c_{\a^{0}\b^{0}}^{\m^{0}} c_{\b^{0}(\a^{1})'}^{\l^{0}}$
is in fact equal to $c_{\m^{0}(\a^{1})'}^{\l^{0}}$.
Similarly, $|\a^n|=0$ and 
$c_{\a^{n-1}\b^{n-1}}^{\m^{n-1}} c_{\b^{n-1}(\a^{n})'}^{\l^{n-1}}$
reduces to
$c_{\a^{n-1}\l^{n-1}}^{\m^{n-1}}$.

\bigskip 
\noindent
\proof
Using (\ref{EQ3}) we obtain
\[
\eta_\m(v) = \sum_{\a^1,\ldots,\a^{n-1}} (-v)^{\sum_i|\a^{i}|}
s_{\m^{0}}(A_0)
\prod_{1\le j\le n-1} s_{\m^{j}/\a^j}(A_j)
s_{(\a^j)'}(A_{j-1})
\,.
\]
Therefore,
\begin{eqnarray*}
\< s_{\ul}\,,\,\eta_\m(v) \>
&=&\sum_{\a^1,\ldots,\a^{n-1}} (-v)^{\sum_i|\a^{i}|}
\<s_{\m^{0}}s_{(\a^1)'}\,,\,s_{\l^{0}}\>
\<s_{\m^{n-1}/\a^{n-1}}\,,\,s_{\l^{n-1}}\> 
\\
&&
\qquad\qquad\qquad\times
\prod_{1\le j\le n-2} 
\<s_{\m^{j}/\a^j}s_{(\a^{j+1})'}\,,\,s_{\l^{j}}\>\,,
\end{eqnarray*}
and the result follows from the identity
$
\<s_{\a/\b}\,s_{\ga}\,,\,s_\de\> 
= \sum_\veps c^{\a}_{\b\veps} c^\de_{\veps\ga}\,
$
and from the fact that $c_{\a\b}^\ga$ is non zero only if 
$|\a|+|\b|=|\ga|$.
\cqfd

\subsection{}
For $\ul=(\l^0,\ldots ,\l^{n-1})\in\P^n$ we set
$\ul'=((\l^{n-1})',\ldots ,(\l^0)')$.
Let $\inf$ denote any total ordering of $\P$ such
that $|\l|<|\m|$ implies $\l\inf\m$. We also denote by $\inf$
the corresponding reverse lexicographic ordering of $\P^n$,
that is, $\ul \inf \um$
if and only if there exists $k$ such that $\l^k\inf\m^k$
and $\l^j=\m^j$ for $j>k$.

The functions $\eta_{\ul}$ and $\vp_{\ul}$ 
enjoy the following properties.

\begin{proposition}\label{prop1}
{\rm (i)} For $\ul,\um \in \P^n$ there holds
$
\< s_{\ul}\,,\,\eta_{\um}(v)  \> =
\<s_{\um'}\,,\, \eta_{\ul'}(v) \>
$.

\smallskip\noindent
{\rm (ii)} For $\ul,\um \in \P^n$ there holds
$
\<\vp_{\ul}(v) \,,\, \eta_{\um}(v)\> = \delta_{\ul\,\um}
$.

\smallskip\noindent
{\rm (iii)}
The matrices
\[
\left[\< s_{\ul} \,,\,\eta_{\um}(v) \>
\right]_{\ul,\um\in\P^n},
\qquad
\left[\< \vp_{\ul}(v) \,,\,s_{\um} \>
\right]_{\ul,\um\in\P^n}
%\qquad
%\left[\langle \vp_{\ul}(v^{-1}) \,,\eta_{\um}(v) \rangle
%\right]_{\ul,\um\in\P^n}
\]
are lower unitriangular if their rows and columns are arranged
according to the total ordering $\sup$.
\end{proposition}
\proof
(i) follows from Lemma~\ref{LEM4} and the symmetry property 
$c_{\a\b}^\ga = c_{\a'\b'}^{\ga'}$ of the Littlewood-Richardson
coefficients.

To prove (ii) we use Lemma~\ref{LEM3}. We have to show that
\[
\sum_{\ul\in\P} \eta_{\ul}(A_0,\ldots ,A_{n-1},v)\,
\vp_{\ul}(B_0,\ldots ,B_{n-1},v)
=
H(A_0B_0)\cdots H(A_{n-1}B_{n-1})\,.
\] 
The left-hand side is equal to
\begin{eqnarray*}
&&\sum_{\l^{0}}s_{\l^{0}}(A_0)
s_{\l^{0}}(B_0+\cdots +v^{n-1}B_{n-1})
\sum_{\l^{1}}s_{\l^{1}}(A_1-vA_0)
s_{\l^{1}}(B_1+\cdots +v^{n-2}B_{n-1})
\\
&&\qquad\quad\times \cdots \times
\sum_{\l^{n-1}}s_{\l^{n-1}}(A_{n-1}-vA_{n-2})
s_{\l^{n-1}}(B_{n-1})
\\
&&
=\quad H(A_0(B_0+\cdots +v^{n-1}B_{n-1}))
H((A_1-vA_0)(B_1+\cdots +v^{n-2}B_{n-1}))
\\
&&
\qquad\quad\times\cdots\times
H((A_{n-1}-vA_{n-2})B_{n-1})
\\
&&
=\quad H(A_0B_0 + \cdots +A_{n-1}B_{n-1})\,.
\end{eqnarray*}
Here we have used (\ref{CAUCHY}) and the basic identities 
\[
H(X)H(Y) = H(X+Y), \quad H(X)H(-X) = 1
\]
valid for any formal sets of variables $X$ and $Y$.  
Hence (ii) is proved.

Suppose that $\<s_{\ul}\,,\,\eta_{\um}(v)\> \not = 0$.
Then, by Lemma~\ref{LEM4} we see that there exists $\a^{n-1}\in\P$
such that 
\[
c_{\a^{n-1}\l^{n-1}}^{\m^{n-1}} \not = 0,
\]
hence either $\l^{n-1}=\m^{n-1}$ or $|\l^{n-1}|<|\m^{n-1}|$.
In the second case $\ul \inf \um$, and in the first one we have
$\a^{n-1} = 0$, therefore
\[ 
c_{\a^{n-2}\b^{n-2}}^{\m^{n-2}}
c_{\b^{n-2}(\a^{n-1})'}^{\l^{n-2}}
= c_{\a^{n-2}\l^{n-2}}^{\m^{n-2}}\not = 0.
\]
It follows that either $\l^{n-2}=\m^{n-2}$ or $|\l^{n-2}|<|\m^{n-2}|$.
In the second case, we obtain that $\ul \inf \um$, and in the first one
we have $\a^{n-2} = 0$.
Thus, by induction we see that 
$\<s_{\ul}\,,\,\eta_{\um}(v)\> \not = 0$.
implies that $\ul \infe \um$. Moreover, clearly
$\<s_{\um}\,,\,\eta_{\um}(v)\> = 1$.
Next, by (ii) we have 
\[
s_{\ul} = \sum_{\um} \<s_{\ul}\,,\,\eta_{\um}(v)\>\, \vp_{\um}(v).
\]
Hence the expansion of $s_{\ul}$ on $\{\vp_{\um}(v)\}$
involves only multi-partitions $\um \infe \ul$, and by solving
this unitriangular system we obtain that
$\<s_{\ul}\,,\,\vp_{\um}(v)\> \not = 0$
implies that $\ul \supe \um$ and
$\<\vp_{\um}(v)\,,\,s_{\um}\> = 1$.
\cqfd

\subsection{}
Define
\[
S=\left[\<s_{\ul} \,,\, \eta_{\um}(v)
  \rangle\right]_{\ul,\um\in\P^n},\qquad
A=S\overline{S^{-1}}=\left[a_{\ul\,\um}(v)\right]_{\ul,\um\in\P^n}.
\]
We introduce a semi-linear involution $f \mapsto \overline{f}$ on $\SS$ by 
requiring that
\[
\overline{s_{\um}} = \sum_{\ul} a_{\ul\,\um}(v)\,s_{\ul}\,,
\qquad
\overline{q(v)\,x} = q(v^{-1})\,\overline{x} \qquad 
(\um\in\P^n, \ x\in\SS, \ q(v) \in K)\,.
\]
The identity $A\overline{S} = S$
shows that 
\begin{equation}\label{equa1}
\overline{\eta_{\um}(v)} = \eta_{\um}(v).
\end{equation}
On the other hand, let $\LL$ be the 
$\Z[v]$-lattice in $\SS$ 
spanned by the vectors $s_{\ul}$.
It follows easily from Lemma~\ref{LEM4} that 
\begin{equation}\label{equa2}
\eta_{\ul}(v) \equiv s_{\ul} \mbox{\ \ mod \ } v\LL\,,\qquad 
\end{equation}   
By Proposition~\ref{prop1} (iii) $S$ is unitriangular,
hence $A$ is also unitriangular. 
So by a classical argument (see \cite{Lu}, 7.10), 
the basis $\{\eta_{\ul}(v)\}$ is uniquely determined
by (\ref{equa1}) and (\ref{equa2}).
This is the canonical basis of $\SS$ associated with the
involution $f \mapsto \overline{f}$ and the lattice $\LL$.
 
\subsection{}
Define
\[
\psi_{\ul}(v^{-1}) = s_{\l^{0}}(A_0)
s_{\l^{1}}(v^{-1}A_0+A_1)
\cdots
s_{\l^{n-1}}(v^{-n+1}A_0+ \cdots +v^{-1}A_{n-2}+ A_{n-1})\,.
\]
Obviously, $\{\psi_{\ul}(v^{-1})\}$ is another basis of $\SS$.
Let $\LL^-$ be the $\Z[v^{-1}]$-lattice in $\SS$ 
spanned by the vectors $s_{\ul}$.
One can check easily that
\begin{equation}\label{equa3}
\psi_{\ul}(v^{-1}) \equiv s_{\ul} \mbox{\ \ mod \ } v\LL^-\,. 
\end{equation}
On the other hand, for $k=1,\ldots,n$ we have the formal identities
\[
v^{-k+1}A_0+ \cdots +v^{-1}A_{k-2}+ A_{k-1}
=[k]A_0 + \sum_{1\le j\le k-1} [k-j] (A_j-vA_{j-1})\,,
\]
where the coefficients $[j]:=(v^j-v^{-j})/(v-v^{-1})$ are
bar-invariant.
Therefore the expansion of $\psi_{\ul}(v^{-1})$ on the
basis $\{\eta_{\um}\}$ will have only bar-invariant coefficients.
Hence
\begin{equation}\label{equa4}
\overline{\psi_{\ul}(v^{-1})} = \psi_{\ul}(v^{-1}),
\end{equation}
and $\{\psi_{\ul}(v^{-1})\}$ is the canonical basis of $\SS$
associated with $f \mapsto \overline{f}$ and the lattice $\LL^-$.

\subsection{}
The coefficients of the expansion of $\psi_{\ul}(v^{-1})$ 
on $\{s_{\ul}\}$ are given by the next lemma, which is proved 
in the same way as Lemma~\ref{LEM4}.
\begin{lemma}\label{LEM6}
For $\ul, \um \in \P^n$ we have
\[
\<\psi_{\ul}(v^{-1})\,,\,s_{\um}\> =
v^{-\De(\ul,\um)}
\sum_{\a} \prod_{0\le k\le n-1}
c^{\m^k}_{\a_k^k,\a_k^{k+1},\ldots , \a_k^{n-1}}
c^{\l^k}_{\a_0^k,\a_1^k,\ldots ,\a_k^k}\,,
\]
where the sum is over all families of partitions
$\a = (\a_i^j\in\P,\  0\le i \le j \le n-1)$
subject to the conditions
\[
\sum_j |\a_i^j| = |\m^i|\,,\qquad \sum_i |\a_i^j| = |\l^j|,
\]
and 
\[
\De(\ul,\um) = \sum_{0\le j \le n-1} j(|\l^j|-|\m^j|)\,.
\]
\end{lemma}

Let $x\mapsto x'$ denote the semi-linear involution of $\SS$
defined by $(s_{\ul})' = s_{\ul'}$.
The next proposition is easily checked by expansion on $\{s_{\ul}\}$. 
\begin{proposition}\label{PROP7}
The adjoint basis $\{\vp_{\ul}(v)\}$ of the canonical basis
$\{\eta_{\ul}(v)\}$ is related to the basis $\{\psi_{\ul}(v^{-1})\}$
by
\[
(\vp_{\ul}(v))' = \psi_{\ul'}(v^{-1})\,. 
\]
\end{proposition}

%%%%%%%%%%%%%%%%%%%%%%%%%%%%%%%%%%%%%%%%%%%%%%%%%%%%%%%%%%%%%
%    SECTION 6
%%%%%%%%%%%%%%%%%%%%%%%%%%%%%%%%%%%%%%%%%%%%%%%%%%%%%%%%%%%%%

\section{Comparison of the canonical bases of $\F_v$ and $\SS$} \label{SECT6}
\subsection{} We fix $w\ge 1$.
Let $\rho = (\rho_1,\ldots ,\rho_l)$ be the large 
$n$-core associated with $w$ (see Definition~\ref{CKcore}).
Let $\P(\rho)$ be the set of partitions with $n$-core $\rho$
and let $\P(\rho,w)$ be the subset of partitions with $n$-weight 
$\le w$.
Note that $\rho$ is equal to its conjugate partition $\rho'$.
It follows that the map $\l \mapsto \l'$ induces a bijection
of $\P(\rho,w)$.

To $\lambda \in \P(\rho)$ we associate its $n$-quotient
$\ul=(\lambda^{0},\ldots ,\lambda^{n-1})\in\P^n$ and its $n$-sign 
$\varepsilon_n(\lambda)$ (see \cite{JK}).
The $n$-quotient of a partition is well-defined only up to circular 
permutations.
We remove this ambiguity by requiring that the partition 
$\rho + (n) = (\rho_1+n,\rho_2,\ldots ,\rho_l)$
have $n$-quotient $(\emptyset,\ldots ,\emptyset, (1))$.
It is well-known that 
$\l \mapsto \ul$
is a bijection from $\P(\rho)$ onto $\P^n$.

Comparing Lemma~\ref{LEM4} with \cite{Mi} we can see that 
for $\l,\m \in \P(\rho,w)$ 
\begin{equation}
\rad_{\l,\m}(v) = \veps_n(\l)\veps_n(\m)\,\< \eta_{\um}(v) \,,\th_{\ul} \>\,
\end{equation}
where $\rad_{\l,\m}(v)$ has been defined in Section~\ref{SECT1}.
This motivates the next definition. 

\subsection{}
Let $\F_v(\rho)$ denote the subspace of $\F_v$ spanned by all $s(\l)$
with $\l\in\P(\rho)$. 
We define a linear map $\Phi$ from $\F_v(\rho)$ to $\SS$ by setting
\[
\Phi(s(\l)) = \veps_n(\l)\,s_{\ul}\,.
\]
This is an isomorphism of vector spaces, which has already been
considered in \cite{LL} for $n=2$ and in \cite{Lei} in general.
Our main result is
\begin{theorem}\label{TH1}
For $\lambda \in \P(\rho,w)$ one has
\begin{eqnarray}
G(\l) &=& \veps_n(\l) \,\Phi^{-1}(\eta_{\ul}(v)) \label{conj1} \,,\\
G^-(\l) &=& \veps_n(\l) \, \Phi^{-1}(\psi_{\ul}(v^{-1}))\,.\label{conj1bis}  
\end{eqnarray}
\end{theorem}
Theorem~\ref{TH1} will be proved in Section~\ref{SECT7}.
\begin{example}
{\rm Let $n=3$, $w=3$, $\rho = (6,4,2,2,1,1)$ and 
$\l =  (12,4,4,3,1,1)$.
Then 
\begin{eqnarray*}
G(\l) &=& s_{(12,4^2,3,1^2)} + v\,s_{(12,4,2^2,1^5)}
+v\,s_{(9,6,5,3,1^2)} + v\,s_{(9,4^2,3^2,2)}\\
&&+\ v^2\,s_{(9,4^2,3,1^5)} + v^2\,s_{(6^2,5,3^2,2)}
+v^2\,s_{(6,4^2,3^2,2^2,1)} + v^3\,s_{(6,4^2,3^2,2,1^3)},
\end{eqnarray*}
therefore
\begin{eqnarray*}
\Phi(G(\l)) & = & -s_{(1)}(A_1)s_{(2)}(A_2) + v\,s_{(1)}(A_0)s_{(2)}(A_2)
+v\,s_{(2)}(A_1)s_{(1)}(A_2) \\
&&+\  v\,s_{(1^2)}(A_1)s_{(1)}(A_2) 
- v^2\,s_{(1)}(A_0)s_{(1)}(A_1)s_{(1)}(A_2)\\
&&-\ v^2\,s_{(2,1)}(A_1) - v^2\,s_{(1^3)}(A_1)
+v^3\,s_{(1)}(A_0)s_{(1^2)}(A_1).
\end{eqnarray*}
On the other hand, 
$(\l^{0},\l^{1},\l^{2}) = (\emptyset, (1),(2))$, $\varepsilon_3(\l)=-1$
and
\begin{eqnarray*}
\eta_{\ul}(v)&=&s_{(1)}(A_1-vA_0)\,s_{(2)}(A_2-vA_1)\\ 
&=&
\left(s_{(1)}(A_1)-v\,s_{(1)}(A_0)
\right)
\left(s_{(2)}(A_2)-v\,s_{(1)}(A_1)s_{(1)}(A_2)
+v^2\,s_{(1^2)}(A_1)
\right)\\
&=& - \Phi(G(\l))\,.
\end{eqnarray*} 
}
\end{example}

\bigskip\noindent
Theorem~\ref{TH1} and Lemma~\ref{LEM4} readily imply
\begin{corollary}\label{formule}
Conjecture~5.5 of {\rm \cite{Mi}} is true and we have for $\l, \m \in \P(\rho,w)$
\[
d_{\l,\m}(v) =
v^{\de(\ul,\um)}
\sum_{\a^{0},\ldots,\a^{n} \atop \b^{0},\ldots,\b^{n-1}}
%c_{\m^{0}(\a^{1})'}^{\l^{0}}
%c_{\a^{n-1}\l^{n-1}}^{\m^{n-1}}
\prod_{0\le j \le n-1}
c_{\a^{j}\b^{j}}^{\m^{j}}
c_{\b^{j}(\a^{j+1})'}^{\l^{j}},
\]
where $\a^{0},\ldots,\a^{n},\b^{0},\ldots,\b^{n-1}$ run through
$\P$ subject to the conditions
\[
|\a^i| = \sum_{0\le j \le i-1} |\l^j|-|\m^j|\,,\qquad
|\b^i| = |\m^i| + \sum_{0\le j \le i-1} |\m^j|-|\l^j|\,,
\] 
and 
\[
\de(\ul,\um) := \sum_{0\le j \le n-2}(n-1-j)(|\l^j|-|\m^j|)\,.
\]
Moroever, we also have
\[
e_{\l,\m}(v) =
v^{\De(\ul,\um)}
\sum_{\a} \prod_{0\le k\le n-1}
c^{\m^k}_{\a_k^k,\a_k^{k+1},\ldots , \a_k^{n-1}}
c^{\l^k}_{\a_0^k,\a_1^k,\ldots ,\a_k^k}\,,
\]
where the sum is over all families of partitions
$\a = (\a_i^j\in\P,\  0\le i \le j \le n-1)$
subject to the conditions
\[
\sum_j |\a_i^j| = |\m^i|\,,\qquad \sum_i |\a_i^j| = |\l^j|,
\]
and 
\[
\De(\ul,\um) = \sum_{0\le j \le n-1} j(|\l^j|-|\m^j|)\,.
\]
\end{corollary}
This can be regarded as a combinatorial description of the
parabolic  Kazhdan-Lusztig polynomials
$d_{\l,\m}(v)$ and $e_{\l,\m}(v)$ in this case.
In particular, we note that these polynomials are just monomials 
when $\l$ and $\m$ belong to $\P(\rho,w)$.

When $n=2$ and $\m$ is a strict partition, $d_{\l,\m}(1)$ is
a decomposition number for the Hecke algebra over a field
of characteristic $0$ at $q=-1$ \cite{LLT96,Ar}. 
In this case, the decomposition numbers for partitions with a large core 
have been first calculated by James and Mathas \cite{JM}.

Finally, we deduce from Proposition~\ref{prop1}~(i) and
Proposition~\ref{PROP7} the following
interesting symmetries
\begin{corollary} For $\l, \m \in \P(\rho,w)$
there holds 
\[
d_{\l,\m}(v) = d_{\m',\l'}(v)\,,\qquad 
e_{\l,\m}(v) = e_{\m',\l'}(v)\,.
\]
\end{corollary}

%%%%%%%%%%%%%%%%%%%%%%%%%%%%%%%%%%%%%%%%%%%%%%%%%%%%%%%%%%%%%
%    SECTION 7
%%%%%%%%%%%%%%%%%%%%%%%%%%%%%%%%%%%%%%%%%%%%%%%%%%%%%%%%%%%%%

\section{Proof of Theorem~\ref{TH1}} \label{SECT7}

\subsection{}
We want to construct certain elements of $U_v(\slchap_n)$ whose
action on $s(\rho) \in \F_v(\rho)$ corresponds via $\Phi$ to the multiplication by 
$e_k(A_j-vA_{j-1})$ in $\SS$ for $k\le w$.
Let $r\in\{0,\ldots , n-1\}$ be the residue of $\rho_1$ modulo $n$.
Thus $F_i s(\rho) = 0$ if $i\not = r$.
For convenience we allow the indices $i$ of the Chevalley generators
$F_i$ to belong to $\Z$ by setting $F_i = F_{i\mod n}$.
\begin{definition}\label{DEFOP}
For $j\in\{0,\ldots ,n-2\}$ and $k\ge 1$ we set
\[
H_{j,k} = F_{j+1+r}^{(k)}\cdots F_{n-2+r}^{(k)}F_{n-1+r}^{(k)}
 F_{j+r}^{(k)}\cdots F_{r+1}^{(k)}F_r^{(k)}\,.
\]
\end{definition}
\begin{example}{\rm Let $n=4$, $w=3$, and $\rho = (12,9,6^2,4^2,2^3,1^3)$.
Then $r=0$ and 
\begin{eqnarray*}
H_{0,k} & = & F^{(k)}_1 F^{(k)}_2 F^{(k)}_3 F^{(k)}_0\,, \\
H_{1,k} & = & F^{(k)}_2 F^{(k)}_3 F^{(k)}_1 F^{(k)}_0\,, \\
H_{2,k} & = & F^{(k)}_3 F^{(k)}_2 F^{(k)}_1 F^{(k)}_0\,. 
\end{eqnarray*}
For $w=2$ we have $\rho = (6,3^2,1^3)$, $r=2$ and
\begin{eqnarray*}
H_{0,k} & = & F^{(k)}_3 F^{(k)}_0 F^{(k)}_1 F^{(k)}_2\,, \\
H_{1,k} & = & F^{(k)}_0 F^{(k)}_1 F^{(k)}_3 F^{(k)}_2\,, \\
H_{2,k} & = & F^{(k)}_1 F^{(k)}_0 F^{(k)}_3 F^{(k)}_2\,. 
\end{eqnarray*}
}
\end{example}
%
%\begin{lemma}\label{Lemma14}
Let $x\mapsto \overline{x}$ denote the bar involution of
$U_v(\slchap_n)$. 
By definition $\overline{F_i^{(k)}} = F_i^{(k)}$ for all $i,k$,
hence
\begin{equation}\label{EqHbar}
\overline{H_{j,k}} = H_{j,k}\,,\qquad (0\le j\le n-2,\ k\ge 1).
\end{equation}
\begin{proposition}\label{prop2}
Let $\l \in \P(\rho)$ with $n$-weight $u<w$.
Then, for $k\le w-u$,
\[
\Phi(H_{j,k} s(\l)) = (-1)^{k(n-j-2)}\veps_n(\l)\, e_k(A_{j+1}-vA_j)s_{\ul}\,.
\]  
\end{proposition}
The proof of Proposition~\ref{prop2} relies on the following combinatorial
properties of the $n$-core partition $\rho$ established in \cite{CK}.
\begin{lemma}\label{LEMCK}
Let $\l,\m\in\P(\rho,w)$ with $|\m|_n = |\l|_n - 1$.
If $\m \subset \l$ 
(\ie the Young diagram of $\m$ is contained in that of $\l$) 
then there exists $i\in\{0,\ldots , n-1\}$ such that $\m^j = \l^j$
for $j\not = i$ and $\m^i \subset \l^i$.
Moreover $\l/\m$ is the Young diagram of the hook partition
$(i+1,1^{n-i-1})$.
\end{lemma}
\begin{figure}[t]
\begin{center}
\leavevmode
\epsfxsize =7cm
\epsffile{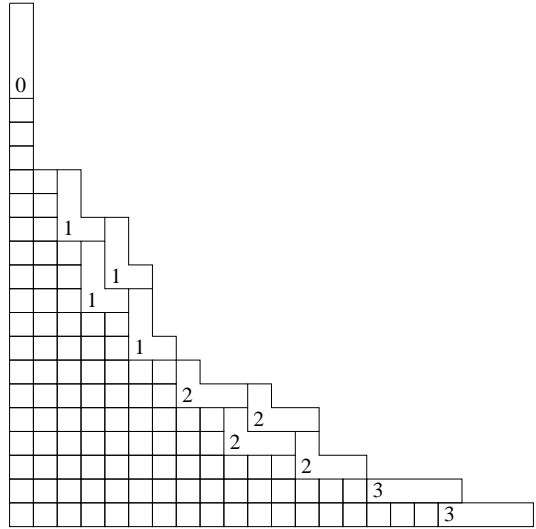}
\end{center}
\caption{\label{Fig3} \small The ribbon tiling of $\l/\rho$ for some
  $\l \in \P(\rho)$. Each ribbon is labelled with its color.}
\end{figure}
\begin{figure}[t]
\begin{center}
\leavevmode
\epsfxsize =13cm
\epsffile{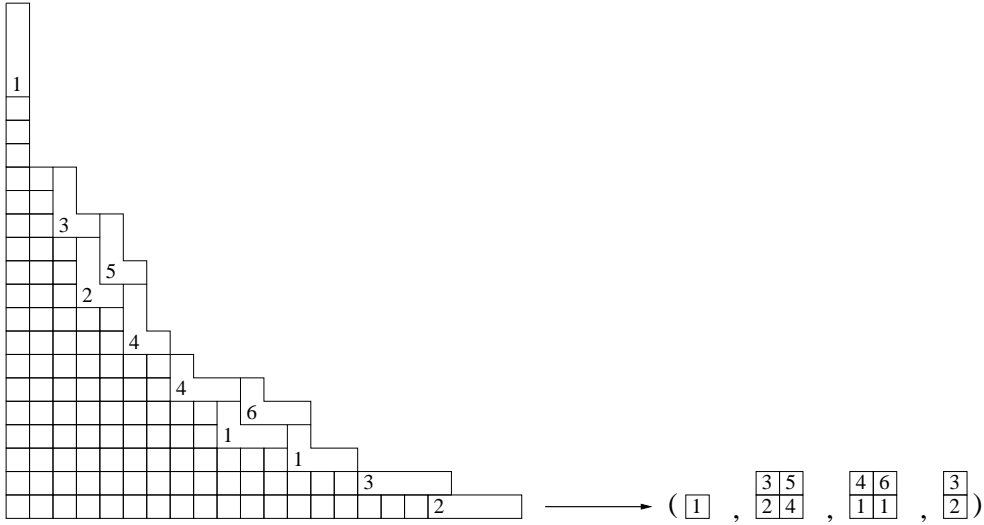}
\end{center}
\caption{\label{Fig4} \small The Stanton-White bijection for a ribbon
tableau of shape $\l \in \P(\rho)$.}
\end{figure}
Thus, if $\l\in\P(\rho,w)$ and if we consider an $n$-ribbon tiling
of $\l/\m$, we see that all ribbons contributing to $\l^i$
have the same shape $(i+1,1^{n-i-1})$ and the same spin $n-i-1$.
Let us say that such ribbons have color $i$.
In fact, for $i<j$ any ribbon of
colour $i$ is situated to the left of any ribbon of colour $j$,
as shown in Figure~\ref{Fig3} with $n=4$. 
Moreover, there is a unique $n$-ribbon tiling of $\l/\m$, and
since the ribbons with different colors do not mix, 
for any ribbon tableau $T$ of shape $\l/\m$ the corresponding
$n$-tuple of tableaux $(t_0,\ldots ,t_{n-1})$ 
under the Stanton-White correspondence \cite{SW,LLT2}
is obtained in a trivial way, as illustrated in Figure~\ref{Fig4}.

\bigskip 
\noindent
{\it Proof of Proposition~\ref{prop2} --- \ }
Since $\wt(H_{j,k}) = -k\de$, all the $s(\m)$ occuring in the
expansion of $H_{j,k}(s(\l))$ have weight $\wt(s(\l))-k\de$, hence
the corresponding partitions $\m$ are obtained from $\l$ by
adding $k$ $n$-ribbons. 
Since by definition
$
H_{j,k} = F_{j+1+r}^{(k)}\cdots F_r^{(k)}
$, 
all these ribbons must be such that one of their ends is
a removable cell of $\m$ with content equal to $j+1+r$ mod~$n$.
With our convention for the definition of the $n$-quotient,
this means that all these ribbons have either colour $j$ or
$j+1$.    
Moreover, the ribbons of colour $j+1$ must form a vertical ribbon strip
and those of colour $j$ must form a horizontal ribbon strip. 
Conversely, let $\m$ be a partition obtained from $\l$ by adding
a horizontal ribbon strip of $s$ ribbons of colour $j$ 
and a vertical ribbon strip of $k-s$ ribbons of colour $j+1$.  
Then there is a unique way of obtaining $\m$ from $\l$ by adding
in that order :
\begin{center}
\begin{tabular}{l} 
$k$ cells with content $r$ mod~$n$, \\ 
$k$ cells with content $r+1$ mod~$n$,\\
$\ldots$\\
$k$ cells with content $r+j$ mod~$n$,\\
$k$ cells with content $n-1+r$ mod~$n$,\\
%$k$ cells with content $n-2+r$ mod~$n$,\\
$\ldots$\\
$k$ cells with content $r+j+1$ mod~$n$.
\end{tabular}
\end{center}
Hence $s(\m)$ occurs in $H_{j,k}(s(\l))$ and its coefficient
is a single power of~$v$.
Then an elementary but tedious calculation based on (\ref{ACTF})
shows that this power is precisely $v^s$. It follows that 
\begin{eqnarray*}
\Phi(H_{j,k} s(\l))& =& \veps_n(\l) \sum_{0\le s \le k}
(-1)^{s(n-1-j)+(k-s)(n-2-j)}
v^s \, h_s(A_j) e_{k-s}(A_{j+1}) \, s_{\ul}\\
&=& (-1)^{k(n-2-j)}\veps_n(\l)\, e_k(A_{j+1}-vA_j)\, s_{\ul}\,,
\end{eqnarray*}
by  \cite{Mcd} (5.16) (5.17).
\cqfd

\subsection{}
Recall the operators $D_k\ (k\ge 1)$ of \ref{q-bosons}. 
\begin{proposition}\label{prop3}
Let $\l \in \P(\rho)$ with $n$-weight $u<w$.
Then, for $k\le w-u$,
\[
\Phi(D_k s(\l)) = 
\veps_n(\l) \, h_k(v^{-n+1}A_0+ \cdots +v^{-1}A_{n-2}+ A_{n-1})\,s_{\ul}\,.
\]   
\end{proposition}
\proof
Let $\m$ be such that $\m/\l$ is a horizontal ribbon strip containing
$k$ ribbons. 
Let $k_i$ be the number of these ribbons which have colour $i$.
Then, by Lemma~\ref{LEMCK}, the spin of $\l/\m$ is equal to
\[
\sp(\l/\m) = \sum_{0\le i \le n-1} k_i(n-1-i)\,,
\]
and the $n$-quotient $\um = (\m^0,\ldots ,\m^{n-1})$ of $\m$ is
obtained from $\ul =(\l^0,\ldots ,\l^{n-1})$ by adding to each
$\l^i$ a horizontal strip of weight $k_i$. 
Using \cite{Mcd} (5.16), it follows that 
\begin{eqnarray*}
\Phi(D_k s(\l))& =& \left(\veps_n(\l) \sum_{k_0+\cdots +k_{n-1}= k}
(-v)^{\sum_i -k_i(n-1-i)}
 \, h_{k_0}(A_0) \cdots h_{k_{n-1}}(A_{n-1})\right) \, s_{\ul}\\
&=&\veps_n(\l) \, \, h_k(v^{-n+1}A_0+ \cdots +v^{-1}A_{n-2}+ A_{n-1})\, s_{\ul}\,.
\end{eqnarray*}
\cqfd

\subsection{} Let $\l\in \P(\rho,w)$ and recall that 
\[
\eta_{\ul}(v) = s_{\l^0}(A_0)s_{\l^1}(A_1-vA_0)\cdots
s_{\l^{n-1}}(A_{n-1}-vA_{n-2})\,.
\]
Each Schur function $s_{\l^i}(A_i-vA_{i-1})$ is a polynomial in
the elementary symmetric functions $e_k(A_i-vA_{i-1})$ with 
coefficients in $\Z$.
On the other hand, using the formal identity
\[
[n]A_0 =
(v^{-n+1}A_0+ \cdots +v^{-1}A_{n-2}+ A_{n-1})
-\sum_{1\le j\le n-1} [n-j] (A_j-vA_{j-1})
\]
where $[j]:=(v^j-v^{-j})/(v-v^{-1})$, we see
that $s_{\l^0}(A_0)$ can be expressed as a polynomial in the variables
$e_k(A_j-vA_{j-1})$ and $h_k(v^{-n+1}A_0+ \cdots +v^{-1}A_{n-2}+
A_{n-1})$
with coefficients in $\C(v)$ invariant under $v\mapsto v^{-1}$.
Therefore by Proposition~\ref{prop2} and
Proposition~\ref{prop3}, since $|\ul|\le w$ the vector 
$\Phi^{-1}\left(\eta_{\ul}(v)\right)$
can be obtained by applying to $s(\rho)$ a polynomial in
the operators $H_{j,k}$ and $D_k$ with bar-invariant coefficients.
Since by (\ref{EqHbar}) and Section~\ref{SECT4.3}
the operators $H_{j,k}$ and $D_k$ are bar-invariant, as well as
$s(\rho)$, it follows that $\Phi^{-1}\left(\eta_{\ul}(v)\right)$
is bar-invariant.
Moreover, we obviously have
\[
\varepsilon_n(\l)\,\Phi^{-1}(\eta_{\ul}(v))
\equiv s(\l) \mod vL\,,
\]
hence (\ref{conj1}) is proved.

Finally, (\ref{conj1bis}) follows easily from (\ref{conj1}).
Indeed, (\ref{conj1}) implies that if $y \in \SS$ is bar-invariant
of degree $\le w$ then $\Phi^{-1}(y) \in \F(\rho)$ is bar-invariant, because
$y$ can be expressed as a linear combination of the 
$\eta_{\ul}(v)$ with bar-invariant coefficients.
Therefore for $\l\in\P(\rho,w)$, 
the vector $\veps_n(\l)\,\Phi^{-1}(\psi_{\ul}(v^{-1}))$ is bar-invariant, and
since it obviously coincides modulo $v^{-1}L^-$ with 
$s(\l)$, it has to be equal to $G^-(\l)$.

Thus, Theorem~\ref{TH1} is proved.  
%%%%%%%%%%%%%%%%%%%%%%%%%%%%%%%%%%%%%%%%%%%%%%%%%%%%%%%%%%%%%%
%    SECTION 8
%%%%%%%%%%%%%%%%%%%%%%%%%%%%%%%%%%%%%%%%%%%%%%%%%%%%%%%%%%%%%

\section{The Scopes isometries} \label{SECT8} 

\subsection{} \label{bijection} 
In \cite{Sc}, Scopes has introduced certain bijections
between sets of partitions with given $n$-cores. 
More precisely, let $\tau$ be an $n$-core partition, and write
$\alpha = \sum_{i=1}^n a_i \veps_i$ for the corresponding element
of $Q_0$, as in \ref{n-cores}. 
Fix $i\in\{0,\ldots ,n-1\}$ and let 
$\si_i(\tau)$ denote the element of $\CC_n$ obtained from
$\tau$ via the action of $W$ described in \ref{n-cores}.
Let $k_i=a_{i+1}-a_i$ if $i\not = 0$ and $k_i=a_1-a_n-1$ if $i=0$.
Thus, $k_i$ is equal to the number of beads transfered
from the $(i+1)$th runner to the $i$th runner when one computes
$\si_i(\tau)$ by means of the abacus representation when $i\not =0$,
and from the 1st to the last runner when $i=0$.
(If $k_i<0$ we understand that the beads are transfered from 
the $i$th runner to the $(i+1)$th runner.)

Let now $\P_0(\tau,w)$ and $\P_0(\si_i(\tau),w)$ be the sets of partitions
with $n$-weight $w$, and $n$-core $\tau$ and $\si_i(\tau)$ respectively.
For $\lambda\in\P_0(\tau,w)$ we deduce from \ref{integrable} and
\ref{scalar} that, $r_i$ being an orthogonal transformation of $\F_\Z$,  
\[
r_i(s(\lambda)) = \pm s(\mu)
\] 
for some $\mu\in\P_0(\si_i(\tau),w)$.
Hence $r_i$ induces a signed bijection $\lambda \mapsto \pm\mu$
from $\P_0(\tau,w)$ to $\P_0(\si_i(\tau),w)$.

\subsection{} \label{w<k}
In the case when $w\le k_i$, it is easy to show that no negative sign occurs
in the previous bijection, as we shall now see.
Let $\sigma(\Lambda_0) \in W\Lambda_0$ be the
extremal weight corresponding to the $n$-core partition~$\tau$,
and let $\Lambda = \sigma(\Lambda_0) - w\delta$.

\begin{lemma} \label{extremal}
%Let $\sigma(\Lambda_0) \in W\Lambda_0$ be the
%extremal weight corresponding to the $n$-core partition~$\tau$,
%and let $\Lambda = \sigma(\Lambda_0) - w\delta$.
Scopes' condition that $w\le k_i$ is equivalent to the fact that
$\Lambda - \alpha_i$ does not belong to~$P(\F)$.
\end{lemma} 
\proof Let us write 
$\sigma(\Lambda_0) = \Lambda_0 + d\delta + \sum_{i=1}^n a_i
\veps_i$.
Then $(\sigma(\Lambda_0)\,,\,\alpha_i) = a_i - a_{i+1}$ if $i\not = 0$,
and $(\sigma(\Lambda_0)\,,\,\alpha_0) = 1+ a_n - a_1$.
Hence, in all cases we have 
$k_i = -(\sigma(\Lambda_0)\,,\,\alpha_i)=-(\Lambda\,,\,\alpha_i)$.
Now
\begin{eqnarray*}
  (\Lambda\,,\,\Lambda)&=&(\sigma(\Lambda_0)-w\delta \,,\,
                               \sigma(\Lambda_0)-w\delta) \\
   &=& (\sigma(\Lambda_0)\,,\,\sigma(\Lambda_0))
       + w^2 (\delta \,,\, \delta)
       -2w (\sigma(\Lambda_0)\,,\,\delta) \\
   &=& (\Lambda_0 \,,\, \Lambda_0) -2w (\Lambda_0 \,,\, \delta) \\
   &=& -2w,   
\end{eqnarray*}
so that $w = -(\Lambda\,,\,\Lambda)/2$.
Therefore, 
\[
w\le k_i 
\quad \Longleftrightarrow \quad 
(\Lambda\,,\,\Lambda) - 2(\Lambda\,,\,\alpha_i) \ge 0
\quad \Longleftrightarrow \quad 
(\Lambda - \alpha_i \,,\, \Lambda - \alpha_i) \ge 1\,.
\]
It is known that $\L\in P(\F)$ if and only if $\L\in \L_0 + Q$
and $(\Lambda\,,\,\Lambda) \le 0$.
Hence $w\le k_i$ is indeed equivalent to $\L-\a_i \not\in P(\F)$.
\cqfd

\begin{lemma}
If $w \le k_i$ , we have
\[
r_i\,x = (e_i^{k_i}/k_i!)\,x\,, \qquad (x\in \F(\L)).
\]
\end{lemma}

\proof
If  $w\le k_i$ by Lemma~\ref{extremal} we have 
$f_i \F(\L) =\{0\}$.
Hence if $x\in \F(\L)$ then 
\[
r_i\,x = \exp(-f_i)\,\exp(e_i)\,\exp(-f_i)\, x
       = \exp(-f_i)\,\exp(e_i) \,x\,.
\]
Since $P(\F)$ is $W$-invariant, Scopes' condition
also implies that $\si_i(\Lambda)+\alpha_i$ does not belong to $P(\F)$.
Therefore we have $e_i \F(\si_i(\Lambda)) = \{0\}$, hence
$e_i^m\,x = 0$ if $m>k_i$. 
Since $r_i\,x \in \F({\si_i(\Lambda)})$ we see that the contributions
of the monomials $e_i^m$ with $m<k_i$ must also cancel, and we 
simply obtain
$r_i\,x = (e_i^{k_i}/k_i!)\,x$.
\cqfd

It follows that if $w\le k_i$, the bijection from $\P_0(\tau,w)$ to
$\P_0(\si_i(\tau),w)$ induced by $r_i$ consists in associating to
$\lambda$ the partition $\mu$ obtained by removing from $\lambda$
the $k_i$ removable nodes with $n$-residue equal to $i$. 
(The previous Lemmas show that this is always possible and in
a unique way.)
So we recover Scopes' description \cite{Sc} of the bijection.
We shall denote this bijection by $\pi_i$.

\subsection{}
Let $A$ be the subring of $K$ consisting of the functions without
pole at $v=0$, and let 
\[
L_A = A\otimes_{\Z[v]}L\,.
\]
By \cite{MM} the $A$-module $L_A$ is a lower crystal lattice at $v=0$ 
in the sense of Kashiwara \cite{Ka1,Ka2}. 
% (this essentially means that $L_A$ is stable under the renormalized
%Chevalley operators $\tilde{E}_i$ and~$\tilde F_i$).
Let $\bar A$ be the subring of $K$ consisting of the functions without
pole at $v=\infty$.
It follows that the $\bar A$-module $\bar L_A$ is a lower
crystal lattice at $v=\infty$.
Finally, let $\F_v^{\inte}$ be the $\C[v,v^{-1}]$-module with basis
$\{G(\lambda)\}$. 
Then we have
\[
\F_v \simeq K\otimes_{A} L_A
\simeq K\otimes_{\,\bar A} \bar L_A
\simeq K\otimes_{\C[v,v^{-1}]} \F_v^{\inte}\,.
\]
Moreover, if we set $E = \F_v^{\inte} \cap L_A \cap \bar L_A$, then clearly
$E$ is a $\C$-vector space with basis $\{G(\lambda)\}$, and 
the map $G(\lambda) \mapsto s(\lambda) \mod vL_A$ is an isomorphism
of $\C$-vector spaces from $E$ to $L_A/vL_A$.
Hence, $(\F_v^{\inte}, L_A, \bar L_A)$ is a balanced triple in the sense
of \cite{Ka3}, and $\{G(\lambda)\}$ is a lower global base of $\F_v$.

In the sequel, in order to simplify notations, we write $\lambda$ instead
of $s(\lambda) \mbox{\ \mod\ } vL_A$. 
The Kashiwara operators $\tilde{E}_i$ and~$\tilde F_i$ act on $L_A/vL_A$.
Since $\{\lambda\}$ is a crystal basis of $L_A/vL_A$ \cite{MM},
for each partition $\lambda$, $\tilde{E}_i\lambda$ 
(\resp $\tilde{F}_i\lambda$) is a single partition $\mu$ or $0$.
The combinatorial description of $\tilde{E}_i\lambda$
(\resp $\tilde{F}_i\lambda$) was given in \cite{MM} (see also \cite{LLT96}).

\subsection{} 
Let $\<\cdot\,,\,\cdot\>_v$ be the scalar product on $\F_v$
defined by 
\[
\< s(\lambda)\,,\,s(\mu)\>_v = v^{-|\lambda|_n}\,\delta_{\lambda,\mu}
\qquad (\lambda, \mu \in \P),
\]
where $|\lambda|_n$ denotes the $n$-weight of $\lambda$.
Let $B^*=\{G^*(\lambda)\}$ be the adjoint basis of $B$
with respect to $\<\cdot\,,\,\cdot\>_v$.
The basis $B^*$ is a renormalization of 
the basis $\{G^\dag(\lambda)\}$ 
of \cite{LT1}, namely 
\[
G^*(\lambda) =
v^{|\lambda|_n}\,G^\dag(\lambda).
\]
Then we have
\begin{equation}\label{vdecomp}
v^{|\lambda|_n} s(\lambda) = \sum_\mu d_{\lambda,\mu}(v) G^*(\mu)\,,
\qquad
G^*(\l) = v^{|\lambda|_n}\,\sum_\m e_{\l',\m'}(-v) s(\m).
\end{equation}

\subsection{} The scalar product $\<\cdot\,,\,\cdot\>_v$ satisfies
\[
\<E_i x\,,\,y\>_v = \<x \,,\,F_i y\>_v, 
\qquad (x,y \in \F_v,\ 0\le i \le n-1),
\] 
(see \cite{LLT96}, 8.1). Therefore $B^*$ is an upper
global base. Hence, using Lemma~5.1.1~(ii) of \cite{Ka3}, we
get
\begin{lemma}\label{maximal}
Let $\lambda \in \P$ and $i\in\{0,\ldots ,n-1\}$.
Let $k$ be the maximal integer such that 
$\tilde{E}^k_i\lambda \not = 0$.
Then $E_i^{(k)} G^*(\lambda) = G^*(\tilde{E}^k_i\lambda)$.
\end{lemma}

\subsection{}
Let us return to the setting of \ref{bijection} and \ref{w<k}.
We fix an $n$-core $\tau$, an integer $i\in\{0,\ldots ,n-1\}$, and
we take $w\le k_i$. 
Then, for all $\lambda \in \P_0(\tau,w)$, the maximal integer $k$ such
that $\tilde{E}_i^k\lambda \not = 0$ is $k=k_i$, and it follows
easily from the combinatorial descriptions of $\tilde E_i \lambda$
and $E_i s(\lambda)$ that
\begin{equation}\label{Eq2}
\tilde{E}_i^{k_i}\lambda = \pi_i(\lambda), 
\quad
E_i^{(k_i)} s(\lambda) = s(\pi_i(\lambda)),
\qquad \lambda \in \P_0(\tau,w),
\end{equation}
where $\pi_i : \P_0(\tau,w) \longrightarrow \P_0(\si_i(\tau),w)$
is Scopes' bijection.
Using Lemma~\ref{maximal} we obtain
\begin{equation}\label{Eq3}
E_i^{(k_i)} G^*(\lambda) = G^*(\pi_i(\lambda)),
\qquad \lambda \in \P_0(\tau,w). 
\end{equation}
Then, combining Equations (\ref{Eq2}) (\ref{Eq3}) (\ref{vdecomp}),
we obtain
\begin{theorem}\label{THSCOPES}
For all $\lambda, \mu \in \P_0(\tau,w)$, there holds
\[
d_{\lambda,\mu}(v) = d_{\pi_i(\lambda),\pi_i(\mu)}(v)\,,\qquad
e_{\lambda,\mu}(v) = e_{\pi_i(\lambda),\pi_i(\mu)}(v)\,,
\]
where $\pi_i$ denotes Scopes' bijection. 
\end{theorem}

\subsection{} Let $\rho=\rho(w)$ and $\si_\rho(\L_0)$
be the large $n$-core associated with $w$ and the corresponding
extremal weight, respectively.
Assume that $\L \in \O_w$ can be reached from $\L_{w,\rho}=\si_\rho(\L_0)-w\de$ by a
sequence of reflections 
\[
\L_{w,\rho} \stackrel{\si_{i_1}}{\longrightarrow}
\L^1  \stackrel{\si_{i_2}}{\longrightarrow} 
\L^2  \stackrel{\si_{i_3}}{\longrightarrow} 
\cdots
 \stackrel{\si_{i_s}}{\longrightarrow} \L^s = \L
\]
such that $\L^j-\a_{i_j} \not \in P(\F)$ for all $j$.
It follows immediately from Theorem~\ref{THSCOPES} that
the transition matrices $T(\L)$ and $T^-(\L)$ are equal to 
$T(\L_{w,\rho})$ and $T^-(\L_{w,\rho})$ respectively and are
given by Corollary~\ref{formule}.

In the case of $\slchap_2$, the orbit $\O_w$ consists of the
weights 
\[
\L^{2k}=(\si_1\si_0)^k(\L_0 - w\de),\qquad
\L^{2k+1}=\si_0(\si_1\si_0)^k(\L_0 - w\de),\qquad
(k\in\N)\,,
\]
and it easy to see that our formulas calculate the canonical bases
for all $\L^j$ with $j\ge w-1$.

In the general case, the class of weights of $\O_w$ to which
the formulas for $\L_{w,\rho}$ can be transferred is also infinite,
but there is still an infinite number of weights of $\O_w$ for which
the formulas do not apply.

\bigskip
\centerline{\bf Acknowledgements}

\bigskip
We want to thank M. Geck and R. Rouquier for organizing a meeting
in Luminy in November 2000 during which our collaboration 
started.
We are grateful to M. Geck, A. Lascoux, R. Rouquier and J.-Y. Thibon 
for stimulating discussions on the problems of this paper.
B.L. acknowledges support from the European Network
{\em Algebraic Lie Representations} ERB FMRX-CT97-0100.

%%%%%%%%%%%%%%%%%%%%%%%%%%%%%%%%%%%%%%%%%%%%%%%%%%%%%%%%%%%%%%%
% BIBLIOGRAPHY
%%%%%%%%%%%%%%%%%%%%%%%%%%%%%%%%%%%%%%%%%%%%%%%%%%%%%%%%%%%%%%%%%
%\newpage

\bigskip
\small

\noindent
\begin{tabular}{ll}
{\sc B. Leclerc} : &
D\'epartement de Math\'ematiques,
Universit\'e de Caen, Campus II,\\
& Bld Mar\'echal Juin,
BP 5186, 14032 Caen cedex, France\\
&email : {\tt leclerc@math.unicaen.fr}\\[5mm]
{\sc H. Miyachi} :&
Department of Mathematics,
Graduate School of Science and Technology,\\
&Chiba University,
Yayoi-cho, Chiba 263-8522, Japan\\
&email : {\tt mmiyachi@g.math.s.chiba-u.ac.jp}
\end{tabular}

\end{document}